\documentclass[amstex,12pt,russian,amssymb]{article}

\usepackage{mathtext}
\usepackage[cp1251]{inputenc}
\usepackage[T2A]{fontenc}
\usepackage[russian]{babel}
\usepackage[dvips]{graphicx}
\usepackage{amsmath}
\usepackage{amssymb}
\usepackage{amsxtra}
\usepackage{latexsym}
\usepackage{ifthen}

\begin{document}

\newcounter{lemma}
\newcommand{\lemma}{\par \refstepcounter{lemma}%
{\bf Лемма \arabic{lemma}.}}

\newcounter{corollary}
\newcommand{\corollary}{\par \refstepcounter{corollary}%
{\bf Следствие \arabic{corollary}.}}

\newcounter{remark}
\newcommand{\remark}{\par \refstepcounter{remark}%
{\bf Замечание \arabic{remark}.}}

\newcounter{theorem}
\newcommand{\theorem}{\par \refstepcounter{theorem}%
{\bf Теорема \arabic{theorem}.}}

\newcounter{proposition}
\newcommand{\proposition}{\par \refstepcounter{proposition}%
{\bf Предложение \arabic{proposition}.}}

\newcommand{\proof}{{\it Доказательство.\,\,}}
\renewcommand{\refname}{\centerline{\bf Список литературы}}

{\bf Е.А.~Севостьянов, А.~Маркиш} (Житомирский государственный
университет им. И.~Франко)

\medskip

{\bf  Є.О.~Севостьянов, А.~Маркиш} (Житомирський державний
університет ім.~І.~Франко)

\medskip

{\bf E.A.~Sevost'yanov, A.~Markysh} (Zhitomir Ivan Franko State
University)

\medskip
{\bf О равностепенной непрерывности одного класса отображений,
квазиконформных в среднем}

\medskip
Изучается поведение открытых дискретных отображений, квазиконформных
в среднем. Доказано, что семейства таких отображений равностепенно
непрерывны (нормальны) в заданной области.

\medskip
{\bf Про одностайну неперервність одного класу відображень,
квазіконформних у середньому }

\medskip
Вивчається поведінка відкритих дискретних відображень,
квазіконформних у середньому. Доведено, що сім'ї таких відображень
одностайно неперервні (нормальні) в заданій області.

\medskip
{\bf On equicontinuity of some class of mappings, which are
quasiconformal in the mean}

\medskip
A behavior of open discrete mappings, which are quasiconformal in
the mean, is investigated. It is proved that the classes of mappings
mentioned above are equicontinuous (normal).

\newpage
{\bf 1. Введение.} В относительно недавних статьях \cite{RS} и
\cite{Sev$_1$} рассмотрены вопросы о локальном поведении классов
отображений, характеристика квазиконформности $Q(x)$ которых
подчинена некоторому интегральному условию, выражаемому при помощи
заданной функцией $\Phi(x)$ (см. там же). Такие отображения мы
называем отображениями, квазиконформными в в среднем. В частности, в
работе \cite{RS} рассмотрены классы гомеоморфизмов, характеристика
которых может меняться (зависеть от отображения), в то время как в
статье \cite{Sev$_1$} <<мажоранта>> $Q(x)$ является фиксированной,
общей для всего рассматриваемого семейства. Основная цель настоящей
заметки -- распространить результаты работы \cite{Sev$_1$} на случай
<<переменной>> мажоранты $Q(x),$ когда класс отображений
определяется только функцией $\Phi(x),$ а не функцией $Q,$
отвечающей за искажение модуля семейств кривых. При этом, усиление
полученных ниже результатов в сравнении с работой \cite{RS} состоит
в том, что здесь рассматриваются открытые дискретные отображения, а
не только гомеоморфизмы.

Основные определения и обозначения, использующиеся ниже, могут быть
найдены в монографии \cite{MRSY} либо статьях \cite{RS} и
\cite{Sev$_1$}.

\medskip
Как обычно, для множеств $E$ и $F\subset\overline{{\Bbb R}^n}$
символ $\Gamma(E,F,D)$ обозначает семейство всех кривых
$\gamma:[a,b]\rightarrow\overline{{\Bbb R}^n},$ которые соединяют
$E$ и $F$ в $D.$ Здесь и далее $h$ -- хордальная метрика (см.
\cite{Va}). Следующая конструкция может быть найдена в работе
\cite{Vu}. Пусть $Q(x, t)=\{y\in \overline{{\Bbb R}^n}: h(x, y)<t\}$
-- сферический шар с центром в точке $x$ радиуса $t.$ Для $x\in
\overline{{\Bbb R}^n},$ множества $E\subset \overline{{\Bbb R}^n}$ и
чисел $0<r<t<1$ полагаем $\widetilde{x}=-\frac{x}{|x|^2},$
$$\left \{\begin{array}{rr}
m_t(E, r, x)= M(\Gamma(\partial Q(x, t), E\cap \overline{Q(x,
r)}))\,,
\\ m(E, x)=m_{\sqrt{3}/2}(E, \frac{\sqrt{2}}{2}, x)\,,
\end{array} \right.$$
и, кроме того,
\begin{equation}\label{eq25}
\left \{\begin{array}{rr} c(E, x)=\max\{m(E, x), m(E,
\widetilde{x})\}\,,
\\ c(E)=\inf\limits_{x\in \overline{{\Bbb R}^n}}c(E, x)\,.
\end{array} \right.
\end{equation}

\medskip
Пусть $\Phi:[0, \infty]\rightarrow [0, \infty]$ -- неубывающая
выпуклая функция. Обозначим через $\frak{R}^{\Phi}_{M,\Delta}$
семейство всех открытых дискретных кольцевых $Q$-отображений в $D,$
таких что $c\left(\overline{{\Bbb R}^n}\setminus f(D)\right)\ge
\Delta$ и
\begin{equation*}\label{eq2!!}
\int\limits_D\Phi\left(Q(x)\right)\frac{dm(x)}{\left(1+|x|^2\right)^n}\
\le\ M\,.
\end{equation*}

\medskip
Имеет место следующая

\medskip
\begin{theorem}\label{th1!}{\it\, Пусть $\Phi:[0, \infty]\rightarrow [0, \infty]$
-- неубывающая выпуклая функция. Если
\begin{equation}\label{eq3!}
\int\limits_{\delta_0}^{\infty}
\frac{d\tau}{\tau\left[\Phi^{-1}(\tau)\right]^{\frac{1}{n-1}}}\ =\
\infty
\end{equation}
для некоторого $\delta_0>\tau_0:=\Phi(0),$ то класс
$\frak{R}^{\Phi}_{M,\Delta}$ является равностепенно непрерывным, и,
следовательно, образует нормальное семейство отображений при всех
$M\in(0, \infty)$ и $\Delta\in(0, 1).$ }
\end{theorem}

\medskip
Ещё раз подчеркнём, что в определении класса
$\frak{R}^{\Phi}_{M,\Delta}$ функция $Q$ не участвует; два различных
отображения $f_1$ и $f_2,$ которым соответствуют разные функции
$Q_1$ и $Q_2,$ могут, вообще говоря, принадлежать одному классу
$\frak{R}^{\Phi}_{M,\Delta}.$

\medskip
{\bf 2. Лемма об оценке искажения.} Пусть $\frak{R}_{Q, \Delta}(D)$
-- класс всех открытых дискретных кольцевых $Q$-отображений $f$ в
области $D\subset{\Bbb R}^n,$ $n\ge 2,$ таких что
$c\left(\overline{{\Bbb R}^n} \backslash f(D)\right) \ge \Delta
>0.$ Для доказательства теоремы
\ref{th1!} установим справедливость следующего утверждения.

\medskip
\begin{lemma}\label{lem1}{\sl\, Пусть $\Delta>0$ и $Q:D \rightarrow [0,\infty]$ --
измеримая функция. Тогда
%for
\begin{equation}\label{eq2*!!}
h\left(f(x),
f(x_0)\right)\le\frac{\omega_{n-1}}{c_n\Delta}\cdot\frac{1}{
\int\limits_{|x-x_0|}^{\varepsilon(x_0)}
\frac{dr}{rq_{x_0}^{\frac{1}{n-1}}(r)}}
\end{equation}
для каждого $f\in \frak{R}_{Q, \Delta}(D)$ и $x\in B(x_0,
\varepsilon(x_0)),$ $\varepsilon(x_0)< {\rm dist} \left(x_0,
\partial D\right),$ где $c_n>0$ зависит только от $n$ и
$q_{x_0}(r)$ -- среднее интегральное значение функции $Q(z)$ на
сфере $|z-x_0|=r$.}
\end{lemma}

\medskip
\begin{proof} Применим подход, использованный
при доказательстве леммы 2 в \cite{Sev$_3$}. Полагаем
$\varepsilon_0:=\varepsilon(x_0).$ Для фиксированного отображения
$f\in {\frak R}_{Q, \Delta}$ рассмотрим конденсаторы $E=(A, C)$ и
$f(E)=E^{\,\prime}=(f(A), f(C)),$ где $C:=\overline{B(x_0,
\varepsilon)},$ $\varepsilon \in(0, \varepsilon_0),$ $A=B(x_0,
\varepsilon_0).$ Условимся, что для конденсатора $E$ символ
$\Gamma_E$ обозначает семейство всех кривых вида
$\gamma:[a,\,b)\rightarrow A,$ таких что $\gamma(a)\in C$ и
$|\gamma|\cap\left(A\setminus F\right)\ne\varnothing$ для
произвольного компакта $F\subset A,$ где $|\gamma|=\{x\in {\Bbb
R}^n: \exists\,t\in [a, b): \gamma(t)=x \}.$ Пусть
$\Gamma_{E^{\,\prime}}$ и $\Gamma_E$ -- соответствующие семейства
кривых для конденсаторов $E$ и $E^{\,\prime},$ соответственно. Ввиду
предложения 10.2 гл. II в \cite{Ri} имеем
$M(\Gamma_{E^{\,\prime}})={\rm cap\,}f(E).$ Полагая
$E_f:=\overline{{\Bbb R}^n}\setminus f(D),$ заметим, кроме того, что
что $\Gamma(f(C), E_f)>\Gamma_{E^{\,\prime}}$ (см. теорему 1.I.46 в
\cite{Ku$_2$}) и, следовательно, в силу свойства минорирования
модуля (см. теорему 6.4 в \cite{Va}) и леммы 1 в \cite{Sev$_2$}
\begin{equation}\label{eq27}
M(\Gamma(f(C), E_f))\le M(\Gamma_{E^{\,\prime}})={\rm cap\,}f(E)\le
\frac{\omega_{n-1}}{I^{n-1}}\,,
\end{equation}
где $\omega_{n-1}$ -- площадь единичной сферы в ${\Bbb R}^n$ и
$$I=I(x_0,\varepsilon,\varepsilon_0)=\int\limits_{\varepsilon}^{\varepsilon_0}
\frac{dr}{rq_{x_0}^{\frac{1}{n-1}}(r)}\,.$$
С другой стороны, согласно теореме 3.14 в \cite{Vu}
\begin{equation}\label{eq28}
M(\Gamma(f(C), E_f))\ge\beta\min\{c(f(C)), c(E_f)\}\,,
\end{equation}
где постоянная $\beta$ зависит только от $n.$ Поскольку для каждого
связного множества $F$ в $\overline{{\Bbb R}^n}$ имеет место
неравенство $c(F)\ge a_n h(F),$ где $h(F)$ -- хордальный диаметр
множества $F,$ а $a_n$ -- некоторая постоянная (см. следствие 3.13 в
\cite{Vu}), будем иметь
\begin{equation}\label{eq29}
c(f(C))\ge a_n \cdot h(f(C))\,.
\end{equation}
Известно, что $c(E)\le
\omega_{n-1}\cdot\left(\log\sqrt{3}\right)^{1-n}$ для любого
множества $E\subset \overline{{\Bbb R}^n}$ (см. соотношение (3.7) в
\cite{Vu} и определение функции $c(\cdot)$ в (\ref{eq25})), так что
\begin{equation}\label{eq30}
\frac{c(E)}{\omega_{n-1}\cdot\left(\log\sqrt{3}\right)^{1-n}}\quad\le\quad
1\qquad \forall\,\, E\subset \overline{{\Bbb R}^n}\,.\end{equation}
Предположим, что $\min$ в правой части (\ref{eq28}) равен $c(f(C)),$
тогда в силу соотношений (\ref{eq29}) и (\ref{eq30})
\begin{equation}\label{eq31}
M(\Gamma(f(C), E_f))\ge\beta\cdot a_n \cdot h(f(C))\ge
\frac{\beta\cdot a_n \cdot
h(f(C))c(E_f)}{\omega_{n-1}\cdot\left(\log\sqrt{3}\right)^{1-n}}\,.
\end{equation}
Пусть $\min\{c(f(C)), c(E_f)\}=c(E_f),$ тогда из (\ref{eq28})
следует, что
\begin{equation}\label{eq32A} M(\Gamma(f(C), E_f))\ge c(E_f)\ge
h(f(C))c(E_f)\,.\end{equation} Полагая $c_n:=\min\left\{1,
\frac{\beta\cdot a_n
}{\omega_{n-1}\cdot\left(\log\sqrt{3}\right)^{1-n}}\right\},$ из
(\ref{eq31}) и (\ref{eq32A}) будем иметь, что
 \begin{equation}\label{eq34}
M(\Gamma(f(C), E_f))\ge c_n\cdot h(f(C))c(E_f)\ge c_n\Delta
h(f(C))\,.
\end{equation}
Из соотношений (\ref{eq27}) и (\ref{eq34}) вытекает, что
\begin{equation}\label{eq35}
h(f(C))\le \frac{\omega_{n-1}}{c_n\Delta I^{n-1}}\,.
\end{equation}
Выберем теперь произвольно $x\in B(x_0, \varepsilon_0),$ $x\ne x_0,$
тогда найдётся $\varepsilon>0,$ $\varepsilon<\varepsilon_0,$ такое
что $|x-x_0|=\varepsilon.$ Поскольку в сделанных выше обозначениях
$x\in C,$ то согласно оценке (\ref{eq35}) мы получим, что
\begin{equation*}\label{eq36}
h(f(x), f(x_0))\le \frac{\omega_{n-1}}{c_n\Delta I^{n-1}}\,,
\end{equation*}
что и доказывает лемму.~$\Box$
\end{proof}

\medskip
{\bf 3. Доказательство теоремы \ref{th1!}.} Идея доказательства
соответствует подходу, использованному при установлении теоремы 4.1
в \cite{RS}. По лемме 3.1 в \cite{RS} имеем оценку
\begin{equation}\label{eq5}
\int\limits_{|x-x_0|}^{\rho}\frac{dr}{rq_{x_0}^{\frac{1}{n-1}}(r)}
=\int\limits_{\varepsilon}^1\frac{dr}{rq^{\frac{1}{n-1}}(r)}\ge
\frac{1}{n}\int\limits_{eM(\varepsilon)}^{\frac{M(\varepsilon)}{\varepsilon^n}}\frac{d\tau}
{\tau\left[\Phi^{-1}(\tau)\right]^{\frac{1}{n-1}}}\,,
\end{equation}
где $\varepsilon=|x-x_0|/\rho,$ $q(r)=q_{x_0}(\rho r)$ и
$$M(\varepsilon)=
\frac{1}{\Omega_n\rho^n\left(1-\varepsilon^n\right)}\int\limits_R
\Phi\left(Q(z)\right)dm(z)\,,$$
$R=\left\{z\in {\Bbb R}^n: |x-x_0|<|z-x_0|<\rho\right\}$ -- кольцо с
центром в точке $x_0$ и $\Omega_n$ -- объём единичного шара ${\Bbb
B}^n$ в ${\Bbb R}^n.$ Т.к. $|z|\le |z-x_0|+|x_0|\le
\rho(x_0)+|x_0|,$ получаем, что
$$M(\varepsilon)\le \frac{\beta_n(x_0)}{\Omega_n
(1-\varepsilon^n)}\int\limits_R
\Phi(Q(z))\frac{dm(z)}{\left(1+|z|^2\right)^n}\,,$$
где $\beta_n(x_0)=\left(1+(\rho(x_0)+|x_0|)^2\right)^n
/\rho^n(x_0).$ Следовательно, при $\varepsilon\le 1/\sqrt[n]{2}$ и,
в частности, при $\varepsilon\le 1/2,$
$$\Phi(0)\le M(\varepsilon)\le \frac{2\beta_n(x_0)}{\Omega_n}M\,.$$
Таким образом, из оценки (\ref{eq5}) вытекает, что для всех $x,$
таких что $|x-x_0|<\rho(x_0)/2,$ имеет место неравенство
\begin{equation}\label{eq6}
\int\limits_{|x-x_0|}^{\rho}\frac{dr}{rq_{x_0}^{\frac{1}{n-1}}(r)}
\ge
\frac{1}{n}\int\limits_{\lambda_n\beta_n(x_0)M}^{\frac{\Phi(0)\rho^n}{|x-x_0|^n}}\frac{d\tau}
{\tau\left[\Phi^{-1}(\tau)\right]^{\frac{1}{n-1}}}\,,
\end{equation}
где $\lambda_n$ -- некоторая постоянная, зависящая только от $n.$
Тогда из (\ref{eq2*!!}) и (\ref{eq6}) вытекает, что
\begin{equation*}\label{eq3}
h\left(f(x),
f(x_0)\right)\le\frac{\omega_{n-1}}{c_n\Delta}\cdot\frac{n}{
\int\limits_{\lambda_n\beta_n(x_0)M}^{\frac{\Phi(0)\rho^n}{|x-x_0|^n}}\frac{d\tau}
{\tau\left[\Phi^{-1}(\tau)\right]^{\frac{1}{n-1}}}}\,\,,
\end{equation*}
откуда ввиду условия (\ref{eq3!}) вытекает равностепенная
непрерывность класса $\frak{R}^{\Phi}_{M,\Delta}$ в точке
$x_0.$~$\Box$

\medskip
КОНТАКТНАЯ ИНФОРМАЦИЯ

\medskip
\noindent{{\bf Евгений Александрович Севостьянов, \\ Антонина Александровна Маркиш} \\
Житомирский государственный университет им.\ И.~Франко\\
ул. Большая Бердичевская, 40 \\
г.~Житомир, Украина, 10 008 \\ e-mail: esevostyanov2009@mail.ru,
tonya@bible.com.ua}

\end{document}